\documentclass{article}
\usepackage[utf8]{inputenc}

\usepackage{titlesec}
\usepackage{amsfonts}
\usepackage{epsfig,graphics,subfigure,psfrag,amsmath,amssymb}
\usepackage{latexsym,wrapfig,picinpar,picins}
\usepackage{bm}
\usepackage{amsthm}
\usepackage{graphicx}
\usepackage{mathrsfs}
\usepackage{fancyref}
\usepackage{dcolumn}
\usepackage{color}
\usepackage{pstricks}
\usepackage{epstopdf}
\usepackage{subfigure}
\usepackage{multicol}
\usepackage{enumerate}
\usepackage{caption}
\usepackage{geometry}
\usepackage{threeparttable}
\usepackage{stmaryrd}
\usepackage{appendix}
\usepackage{cases}
\usepackage[section]{placeins}
 \usepackage{ulem}

\makeatletter
\@addtoreset{equation}{section}
\makeatother

\newcommand\eps{\varepsilon}
\newcommand\be{\begin{equation}}
\newcommand\ee{\end{equation}}
\newcommand\ben{\begin{equation*}}
\newcommand\een{\end{equation*}}
\newcommand\bea{\begin{eqnarray}}
\newcommand\eea{\end{eqnarray}}
\newcommand\bean{\begin{eqnarray*}}
\newcommand\eean{\end{eqnarray*}}
\newcommand\pa{\partial}

\newcommand\no{\nonumber \\}

\titleformat*{\section}{\Large\bfseries}
\titleformat*{\subsection}{\normalsize\bfseries}
\begin{document}
\newtheorem{remark}{\emph Remark}[section]
\newtheorem{lemma}{\emph Lemma}[section]
\newtheorem{definition}{\emph Definition}[section]
\newtheorem{theorem}{\emph Theorem}[section]

\title{Homogenization Theory of Ion Transportation in Multicellular Tissue}

\author{Chun Xiao \footnotemark[1]  \and Xingye Yue\footnotemark[3]  \and
Huaxiong Huang\footnotemark[4]\and Shixin Xu \footnotemark[2] 
}

 \renewcommand{\thefootnote}{\fnsymbol{footnote}}
 \footnotetext[1]{School of Mathematical Sciences, Soochow University, Suzhou 215006, China; School of Mathematics and Statistics, Lingnan Normal University, Zhanjiang 524048, China (cxiao77xc@gmail.com).}
 \footnotetext[2]{Corresponding author. Zu Chongzhi Center for Mathematics and Computational Sciences (CMCS), Duke Kunshan University, Kunshan 215316, China (shixin.xu@dukekunshan.edu.cn).}
 \footnotetext[3]{Department of Mathematics, Soochow University, Suzhou 215006, China (xyyue@suda.edu.cn).}
 \footnotetext[4]{Research Center for Mathematics, Advanced Institute of Natural Sciences, Beijing
Normal University (Zhuhai), Guangdong, China;
Department of Mathematics and Statistics, York University, Toronto, ON, Canada (hhuang@uic.edu.cn).}


\date{}
\maketitle

{
\begin{abstract}
\footnotesize
  Ion transport in biological tissues is crucial in the study of many biological and pathological problems. Some multi-cellular structures, like smooth muscles on the vessel walls, could be treated as periodic bi-domain structures, which consist of intracellular space and extracellular space with semipermeable membranes in between. With the aid of two-scale homogenization theory, macro-scale models
are proposed based on an electro-neutral (EN) microscale model with nonlinear interface conditions,
where membranes are treated as combinations of capacitors and resistors. The connectivity of intracellular space is also taken into consideration. If the intracellular space is fully connected and forms a
syncytium, then the macroscale model is a bidomain nonlinear coupled partial differential equations system. Otherwise, when the intracellular cells are not connected, the macroscale model for intracellular space is an ordinary differential system with source/sink terms from the connected extracellular space.


 \end{abstract}
}

\section{Introduction}
\par Ions in human body play vital roles in many aspects such as helping the transport of nerve impulses, maintaining the proper functions of muscles, activating various enzymes, helping blood coagulation and so on. Studying ion transport in biological tissues can help us understand the mechanisms of many physiological phenomenon and gain some insight about how to treat certain diseases. The Poisson-Nernst-Planck (PNP) system is one of the most popular  mathematical models that describe  the ion transport under the influence of both ionic concentration gradient and electric field. PNP system has extensive and successful applications in biological systems, particularly in ion channels on cell membrane
\cite{cardenas2000three,chen1997permeation,corry2000tests,Gillespie(2001),Horng(2012),nonner1998ion}. Due to the capacitance of membranes, there are  thin boundary layers (BLs) near the interfaces formed by excessive charges accumulation.
BLs  requires  extra computation cost  during numerical simulations in order to resolve the fast change behaviors of solution inside the layers and  attain certain accuracy \cite{budd1996moving,tang2003adaptive}. A lot of efforts are put in order to get rid of this constrain, like
Mori \cite{mori2009three}.  By using asymptotic analysis, Song et al. \cite{Song(2018)1, Song(2018)2}  proposed effective interface conditions  by introducing a time dependent capacitance on the membrane.  In this paper, we take the linearization of the effective boundary conditions and propose a microscale EN PNP system with interface conditions describing the membrane fluxes and capacitor effect.

Due to the existence of pumps on membranes, ion concentrations across the membrane are discontinuous, for example potassium is $140\sim 150 ~mM$ in the intracelluar space and $3\sim 5~ mM$ in the extracellular space. However, the flux across the membrane is continuous and determined by conductance and the difference between membrane potential and Nernst Potential. When studying biological tissues, there are thousands of cells in the system and the solution is highly oscillatory.
The obtained mircoscale model requires significant computational resources to solve numerically. In order to simulate ion transport in biological tissues more efficiently, an effective macroscale model is demanded.
 One of the most popular ways to derive the macroscale models is  through homogenization by deriving an ``average'' or ''effective'' homogenized system, which extracts macro information from micro structures. Specific methods of homogenization include oscillating test function method \cite{tartar1976quelques}, asymptotic expansion method \cite{bakhvalov1992averaging}, two-scale method \cite{nguetseng1989general,Allaire(1992)} and unfolding method \cite{Cioranescu(2006),Cioranescu(2002)}.
The homogenization theory for system with jump solution is first developed by  Monsurro  and Donato \cite{monsurro2002homogenization,monsurro2004erratum,donato2004homogenization,donato2006some} to  a linear elliptic model for heat conductors problem by using extension operators.  
Results for linear parabolic and hyperbolic equations can be found in \cite{Jose2009Transport,donato2010corrector,yang2014periodic,yang2014homogenization,donato2009correctors}.
 Another factor for bi-domain homogenization is whether  two sub-domains are all connected domain, respectively. The linear problems mentioned above only considered the case when one sub-domain  is embedded in the other sub-domain and disconnected. In \cite{Bunoiu(2018)}, a linear diffusion problem in a bi-domain with flux jump at the interface is discussed. The authors considered both when one sub-domain is connected and disconnected, and explained the reasons for the difference in the resulting homogenization systems for the two cases, which is mainly due to the different spaces of the test functions. For nonlinear bi-domain homogenization problems, the main difficulty is the strong convergence requirement both in domain and on the interface.
 By two-scale convergence and unfolding operator, Gahn et al. \cite{graf2014homogenization, Gahn(2016)}   proposed homogenization results for 
 a nonlinear problem in a bi-domain  for calcium dynamics (connected sub-domains) and diffusion-reaction system (one disconnected sub-domain).

For the homogenization theory of PNP system, which is a nonlinear coupled system,   asymptotic expansion method \cite{looker2006homogenization} is first used to derive the homogenized system for  the linearized Navier-Stokes-PNP model in porous media in \cite{o1978electrophoretic}. Later, a  rigorous derivation by using two-scale method is proposed by Alliare et al. \cite{allaire2010homogenization}. Similar method is extended to study ion transport through deformable porous media \cite{allaire2017ion}.
Homogenization of a full nonlinear PNP model in porous media is discussed in \cite{timofte2014homogenization,schmuck2015homogenization} by using unfolding method and two-scale method, where nonflux boundary condition is used on the interface for ion concentraion. 
 Most of the homogenization research for PNP model and ion transport model are established in porous media without considering the electric-diffusion of ions in the intracellular region, and sometimes linearization technique is used to simplify the process.   In this paper, by using unfolding operator and two-scale convergence,  we develop the homogenization theory for the fully nonlinear EN bi-domain model  in the whole domain with nonlinear interface  flux boundary conditions which  depends on the jumps of ion concentrations, electric potential,  and  the time derivative of potential jump. 
 Two different scenarios, connected and disconnected intracellular regions, are both taken into consideration and lead to different macroscopic models.  

The remainder of this paper is organized as follows: The microscale EN ion transport model is given in Section \ref{section EN model}; Section \ref{section a prior} is devoted to proving the a priori estimates; Then convergence results and homogenization process are presented in Section \ref{section homo 2} according to different connectivity conditions of intracellular region; we draw conclusions in the last section.

\section{Setting of the mathematical model}  \label{section EN model}
\begin{figure}[ht]
  \centering
  \includegraphics[width=0.7\textwidth]{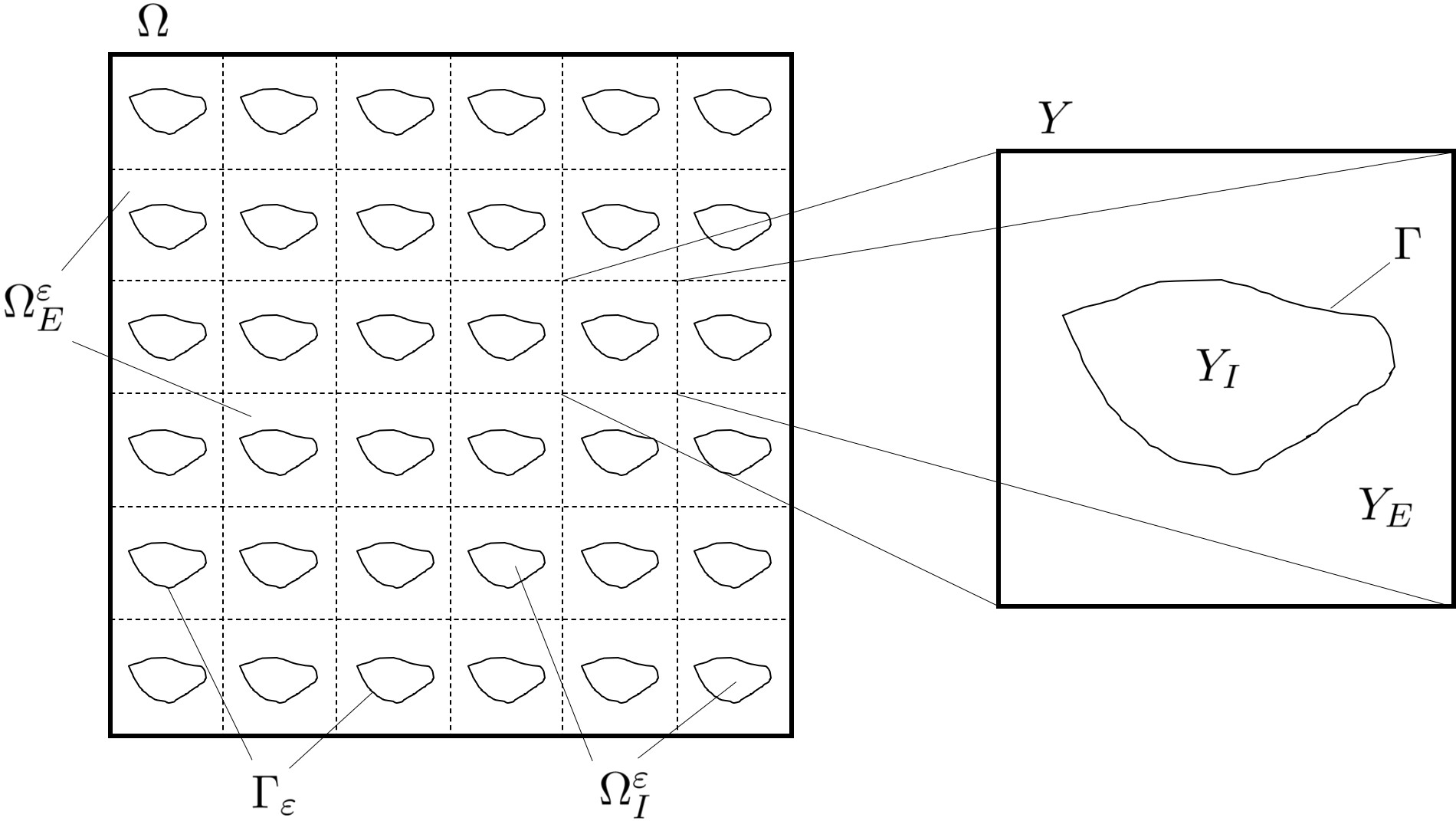}
  \caption{{\footnotesize Schematic of the domain $\Omega$.}}\label{fig tissue01}
\end{figure}

In this section, a microscale EN bi-domain ion transport model in multicellular tissues is introduced. Consider a domain $\Omega=(0,1)^d$ which consists of two components: $\Omega_I^\eps$ and $\Omega_E^\eps$ (see Figure \ref{fig tissue01}). Let $Y=(0,1)^d$ and $Y_I,Y_E$ are two disjoint subsets of $Y$, such that
\ben
\overline{Y}=\overline{Y_{E}} \cup \overline{Y_{I}}.
\een
And $\Gamma=\partial Y_{I} \bigcap \partial Y_{E} $ is Lipschitz continuous. Let $\vec{n}_1$ be the normal direction of $\Gamma$, pointing from $Y_I$ to $Y_E$. Let $\beta_s, s=I,E$ be the characteristic function on $Y_s,s=I,E$, which extend periodically to $\mathbb{R}^d$. For any $k \in \mathbb{Z}^{d}$, let
\ben
Y^{k}=k+Y, ~\Gamma^{k}=k+\Gamma, ~Y_{s}^{k}=k+Y_{s},
\een
where $k=\left(k_{1} , \cdots, k_{d} \right)$, $s=I,E$. For any $\eps>0$ and $1/\eps \in \mathbb{N}^+$, let $K_{\varepsilon}=\left\{k \in \mathbb{Z}^{d} | \varepsilon Y_{s}^{k} \cap \Omega \neq \emptyset, s=I,E\right\}$.
Denote the two disjoint subsets of $\Omega$ and the interface between them as:
\ben
\Omega_{I}^\eps=\bigcup_{k \in K_{\varepsilon}} \varepsilon Y_{I}^{k}, \quad \Omega_{E}^\eps=\Omega \backslash \overline{\Omega_{I}^\eps}, \quad \Gamma_{\varepsilon}=\bigcup_{k \in K_{\varepsilon}} \varepsilon \Gamma^{k}.
\een
So $\Omega=\Omega_I^\eps \bigcup \Omega_E^\eps \bigcup \Gamma_\eps$ and $\Omega$ is a union of $\eps Y$-periodic sets. 
 Assume both $\Omega_E^\eps$ and $\Omega_I^\eps$ have Lipschitz boundaries, especially $\Gamma_\eps$ is a Lipschitz boundary.

Consider three different ion species in $\Omega$: $\mathrm{Na^{+},K^{+},Cl^{-}}$, their concentrations and valences are denoted by $C_{i,s}^\eps, i=1,2,3, s= I,E$ and $z_i$ respectively, the electric potential is $\phi_s^\eps$, subscript $s$ represent variables in $\Omega_s^\eps,s=I,E$. 

The assumptions are summarized as follows:
 \begin{itemize}
 \item Connectivity: $\Omega_E^\eps$ is connected and $\Omega_I^\eps$ can be connected or not (see Figure \ref{fig tissue01}).
\item Positivity:
\bea 
0< C_d \leq C_{i,s}^\eps \leq C_u~~~~\text{on} ~\Gamma_\eps,~~i=1,2,3,s=I,E,\label{assum positive}\\
   0<C_l<\sum\limits_{i=1}^{3} z_i^2 C_{i,s}^\eps (t,x)~~~~\text{in}  ~\Omega_s^\eps,~~s=I,E,\label{assum positive 01}
\eea
where $C_d,C_u,C_l$ are positive constants.
\item Electric Neutrality (EN):
\be   \label{assum elecnuetra}
\sum\limits_{i=1}^{3} z_i C_{i,s}^\eps (t,x)=0,~~s=I,E.
\ee
\item Diffusion Constant: diffusion constants of ions are the same and denoted by  $D$.
 \end{itemize}
 When $\Omega_I^\eps$ is not connected, we call this case ``connected-disconnected'' case. In this case, $\overline{Y_I} \subset Y$ and $\Omega_I^\eps$ is a disconnected union of $\eps Y$-periodic sets of size $\eps Y_I$.
When $\Omega_I^\eps$ is connected, we call this case ``connected-connected'' case. In this case both $Y_I$ and $Y_E$ reach $\pa Y$, thus both $\Omega_I^\eps$ and $\Omega_E^\eps$ reach $\pa \Omega$. 

 The EN bi-domain ion transport model is given as follow: for $i=1, 2,3,s=I,E$,
\begin{subequations} \label{elec neu nondimen two ions}
\begin{align}
\pa_t C_{i,s}^\eps=-\nabla\cdot {J}_{i,s}= \nabla\cdot[D( \nabla C_{i,s}^\eps+ z_i C_{i,s}^\eps \nabla \phi_{s}^\eps)]
~~&in ~ (0,T)\times \Omega^\eps_s  , \label{elec neu nondimen two ions sub a}\\
 z_i {J}_{i,s} \cdot \vec{n}
= \eps G_i\big( \llbracket{\phi^\eps }\rrbracket +
\dfrac{1}{z_i } \llbracket{\ln C_i^\eps}\rrbracket \big)+\eps P_i^\eps+\eps \lambda_i P_m \pa_t \llbracket{\phi^\eps }\rrbracket
 ~~ & on ~ (0,T)\times\Gamma_\eps, \label{elec neu nondimen two ions sub b}\\
{J}_{i,E}\cdot \vec{\nu}=0 ~~& on~ (0,T)\times\pa \Omega,\label{elec neu nondimen two ions sub c}\\
C_{i,s}^\eps|_{t=0}=C_{i,s}^0 ~~& in~\Omega^\eps_s.\label{elec neu nondimen two ions sub d}
\end{align}
\end{subequations}
and
\begin{subequations}  \label{potential nodimen two ion}
\begin{align}
-\nabla \cdot (D\sigma_{s}^\eps  \nabla \phi_{s}^\eps )=0  ~~&   in ~  (0,T)\times\Omega^\eps _s, \label{potential nodimen two ion sub a} \\
-D\sigma_{s}^{\eps} \nabla \phi_s^\eps  \cdot \vec{n}=
\sum\limits_{i=1}^{3}\eps G_i\llbracket{\phi^\eps  }\rrbracket+\sum\limits_{i=1}^{3} \eps G_i\dfrac{1}{z_i}\llbracket{\ln C_i^\eps }\rrbracket+\eps I_P^\eps+\eps P_m \pa_t \llbracket{\phi^\eps  }\rrbracket     ~~&  on ~ (0,T)\times\Gamma_\eps  , \label{potential nodimen two ion sub b}\\
-D\sigma_{E}^{\eps} \nabla \phi_E^\eps  \cdot \vec{\nu}=0   ~~&   on~ (0,T)\times\pa \Omega, \label{potential nodimen two ion sub c}\\
\llbracket{\phi^\eps}\rrbracket|_{t=0}=\phi^0 ~~& on~\Gamma_\eps. \label{potential nodimen two ion sub d}
\end{align}
\end{subequations}
where   $\phi_s^\epsilon$ is the electric potential in domain $\Omega_s$,  $C_{i,s}^\epsilon$ is the concentration of $i_{th}$ ion  in domain $\Omega_s$.

\eqref{elec neu nondimen two ions sub a} is the usual Nernst-Planck equation, \eqref{potential nodimen two ion} is a direct result of \eqref{elec neu nondimen two ions} and EN assumption \eqref{assum elecnuetra} with 
\be \label{sigma_s relation}
\sigma_{s}^\eps=\sum\limits_{i=1}^3 z_i^2 C_{i,s}^\eps,~~s=I,E.
\ee
\eqref{elec neu nondimen two ions sub b} is the  interface flux condition  for ion concentration which consists of three parts.
\begin{itemize}
    \item The first part is the current induced by passive ion channel which is modeled by Ohm's Law. Here $G_i$ is the conductance of $i_{th}$ ion on the membrane, $\llbracket{\phi^\eps }\rrbracket  =\phi_{I}-\phi_E$ is the membrane potential and $E_i \triangleq -\frac{1}{z_i} \llbracket{\ln C_i^\eps }\rrbracket=\frac{1}{z_i} (\ln C_{i,E}^\eps-\ln C_{i,I}^\eps)$ is the Nernst potential. 
    \item The second part is the current induced by the active pumps on the membrane. Here we only consider the $Na^{+}/K^{+}$ pumps \cite{zhu2019bidomain}.
\ben
P_1^\eps=3 I_p^\eps~~,P_2^\eps=-2 I_p^\eps,~~P_3^\eps=0,
\een
with
\begin{eqnarray}  \label{I_p express}
I_p^\eps&=&
I_{\max 1}\left(\frac{C_{1,I}^\eps}{C_{1,I}^\eps+K_{N a 1}}\right)^{3}\left(\frac{C_{2,E}^\eps}{C_{2,E}^\eps+K_{K 1}}\right)^{2}\no
&+&
I_{\max 2}\left(\frac{C_{1,I}^\eps}{C_{1,I}^\eps+K_{N a 2}}\right)^{3}\left(\frac{C_{2,E}^\eps}{C_{2,E}^\eps+K_{K 2}}\right)^{2},
\end{eqnarray}
where $I_{max1},I_{max2}$ are maximum $Na^{+}/K^{+}$ pump current density and $K_{K 1},K_{K 2}$, $K_{Na 1}$, $K_{Na 2}$ are threshold constants for $I_p^\eps$. 
\item The last term is the current induced by the capacitor effect of membrane with capacitance $P_m$   .  $\lambda_i$ reflects the effect of capacitor on the $i_{th}$ ion species and $\sum\limits_{i=1}^3 \lambda_i=1$. Here we take $\lambda _i$ as constant \cite{mori2009three} which is simplification and linearization version of models in  \cite{Song(2018)1, Song(2018)2}. 
\end{itemize}

 \begin{remark}
In \eqref{elec neu nondimen two ions}, \eqref{potential nodimen two ion} we only consider the case when $\Omega_I^\eps$ is not connected, and this would not affect the main idea because when $\Omega_I^\eps$ is also connected and reaches the exterior boundary $\pa \Omega$, we only need to add another exterior boundary condition for the functions in $\Omega_I^\eps$. And the connectivity condition of $\Omega_I^\eps$ would not affect the a priori estimates (Theorem \ref{Th a prior est}) either.
\end{remark}
The homogenization results (Theorem \ref{theorem con-discon} and Theorem \ref{two ions converge th 01})  show  that the interfacial fluxes become source terms in the homogenized equations, and the connectivity of intracellular region will affect the homogenization results.  When the intracellular region is connected, the homogenized one is a diffusion-reaction equation; when cells are isolated with each other, the homogenized equation is a reaction equation. This  is consistent with relevant results in  \cite{Bunoiu(2018)}. And the EN condition still holds for the homogenized equations both when intracellular region is connected and not connected.

\section{A priori estimate}  \label{section a prior}
In order to derive the convergence results, we need a priori estimates for the solutions of \eqref{elec neu nondimen two ions}, \eqref{potential nodimen two ion}. In the rest of the paper, $A$ (with or without subscript) represent a generic positive constant independent of $\eps$ (in different equations, the value of $A$ may be different).

The weak form of \eqref{elec neu nondimen two ions} and \eqref{potential nodimen two ion} is: for $i=1,2,3,s=I,E$, find $C_{i,s}^\eps\in L^2(0,T;~H^1(\Omega_s^\eps)),~\phi_s^\eps \in L^2(0,T;~H^1(\Omega_s^\eps)) $, such that
\begin{subequations}  \label{weak form concentration two ion}
\begin{align}
& \int_{\Omega_s^\eps} \pa_t C_{i,s}^\eps v_s~dx+\int_{\Omega_s^\eps} D \nabla C_{i,s}^\eps \nabla v_s~dx +z_i \int_{\Omega_s^\eps}
C_{i,s}^\eps D \nabla \phi_s^\eps \nabla v_s~dx \no
&\pm \eps  \int_{\Gamma_\eps} \big(\dfrac{G_i}{z_i}\llbracket{\phi^\eps}\rrbracket v_s+\dfrac{G_i}{z_i^2}\llbracket{\ln C_i^\eps}\rrbracket v_s+\dfrac{P_i^\eps}{z_i} v_s+
\dfrac{\lambda_{i}}{z_i}P_m \pa_t\llbracket{\phi^\eps}\rrbracket v_s \big)~dS=0 ,\label{weak form concentration two ion sub a}\\
& \int_{\Omega_s^\eps} \sigma_s^\eps D \nabla \phi_s^\eps \nabla \varphi_s~dx \pm \eps
  \int_{\Gamma_\eps} \big( \sum\limits_{k=1}^{3}G_k\llbracket{\phi^\eps}\rrbracket\varphi_s+\sum\limits_{k=1}^{3}\dfrac{G_k}{z_k}
\llbracket{\ln C_k^\eps }\rrbracket\varphi_s +I_p^\eps \varphi_s +P_m \pa_t\llbracket{\phi^\eps}\rrbracket\varphi_s \big)~dS=0 ,\label{weak form concentration two ion sub b}
\end{align}
\end{subequations}
for $\forall ~v_s \in L^2(0,T;~H^1(\Omega_s^\eps))$ and $\varphi_s \in L^2(0,T;~H^1(\Omega_s^\eps))$, where ``$\pm$'' takes ``$+$'' for $s=I$ and ``$-$'' for $s=E$. Next we derive a priori estimates for $C_{i,s}^\eps,\phi_s^\eps,~i=1,2,3,s=I,E$.

We also impose the following condition for $\phi_E^\eps$ to ensure the uniqueness, 
\be \label{phi_s assump zero mean}
\int_{\Omega_E^\eps} \phi_E^\eps ~dx=0.
\ee

In the rest of the paper, for a subset $\omega \in \mathbb{R}^d$, we simply denote $\|f\|_{L^2(\omega)}$ as $\|f\|_\omega$. And for a function $u^\eps$ defined both in $\Omega_I^\eps$ and $\Omega_E^\eps$, denote $ \| u^\eps \|^2_{\Omega^\eps} \triangleq \| u_{I}^\eps \|^2_{\Omega_I^\eps}+\| u_{E}^\eps \|^2_{\Omega_E^\eps} $.
Then a priori estimates are summarized in the following theorem.
\begin{theorem}  \label{Th a prior est}
Let $C_{i,s}^\eps,\phi_s^\eps,i=1,2,3,s=I,E$ be the solutions of \eqref{weak form concentration two ion} and suppose assumption \eqref{phi_s assump zero mean} hold, then the following estimates are valid, for $i=1,2,3,s=I,E$,

\begin{subequations}  \label{th prior result}
\begin{align}
& \| C_{i,s}^\eps \|_{L^\infty(0,T;L^2(\Omega_s^\eps))}\leq A ,\label{th prior result sub a}\\
 &\| C_{i,s}^\eps \|_{L^2(0,T;H^1(\Omega_s^\eps))} \leq A,\label{th prior result sub b}\\
&\sqrt{\eps} \| C_{i,s}^\eps  \|_{L^2(0,T;L^2(\Gamma^\eps))}  \leq A,\label{th prior result sub c}\\
& \|\pa_t C^\eps_{i,s} \|_{L^2(0,T;H^{-1}(\Omega^\eps_s))}   \leq A,\label{th prior result sub d}\\
&\sqrt{\eps} \| \phi^\eps  \|_{L^2(0,T;L^2(\Gamma^\eps))}   \leq A,\label{th prior result sub e}\\
&\| \phi_s^\eps \|_{L^2(0,T;H^1(\Omega_s^\eps))}  \leq A. \label{th prior result sub f}
\end{align}
\end{subequations}
\end{theorem}
\textit{Proof}: 
 If let $v_s=C_{i,s}^\eps,~\varphi_s=\phi_s^\eps,~s=I,E$ in \eqref{weak form concentration two ion}, we have
\begin{subequations}  \label{prior middle01}
\begin{align}
\hspace{-20pt} \sum\limits_{i=1}^{3} \dfrac{1}{2} \pa_t \| C_i^\eps \|^2_{\Omega^\eps}+
\sum\limits_{i=1}^{3} D \| \nabla C_i^\eps \|^2_{\Omega^\eps} +
\sum\limits_{i=1}^{3} \int_{\Omega^\eps} z_i C_i^\eps D \nabla \phi^\eps \nabla C_i^\eps ~dx +\sum\limits_{i=1}^{3}\int_{\Gamma^\eps} \eps \dfrac{P_i^\eps}{z_i} \llbracket{C_i^\eps }\rrbracket~dS\no
\hspace{-20pt}+\sum\limits_{i=1}^{3}\int_{\Gamma^\eps} \eps \dfrac{G_i}{z_i} \llbracket{\phi^\eps }\rrbracket
 \llbracket{C_i^\eps }\rrbracket~dS+
\sum\limits_{i=1}^{3}\int_{\Gamma^\eps} \eps \dfrac{G_i}{z_i^2} \llbracket{\ln C_i^\eps }\rrbracket
\llbracket{ C_i^\eps }\rrbracket ~dS\no
\hspace{-20pt} +
\sum\limits_{i=1}^{3}\int_{\Gamma^\eps} \eps \dfrac{\lambda_i}{z_i} P_m \pa_t
\llbracket{\phi^\eps }\rrbracket \llbracket{C_i^\eps }\rrbracket~dS
=0 , \label{prior middle01 sub a}\\
\hspace{-20pt}  \int_{\Omega^\eps} \sigma^\eps D \nabla \phi^\eps \nabla \phi^\eps~dx +
\sum\limits_{i=1}^{3}\eps G_i \| \llbracket{\phi^\eps }\rrbracket \|^2_{\Gamma^\eps}+
\eps\dfrac{P_m}{2}\pa_t\| \llbracket{\phi^\eps }\rrbracket \|^2_{\Gamma^\eps}\no
\hspace{-20pt}=-\sum\limits_{i=1}^{3} \int_{\Gamma^\eps} \eps \dfrac{G_i}{z_i} \llbracket{\ln C_i^\eps }\rrbracket \llbracket{\phi^\eps }\rrbracket~dS
-\int_{\Gamma^\eps} \eps I_p^\eps  \llbracket{\phi^\eps }\rrbracket~dS. \label{prior middle01 sub b}
\end{align}
\end{subequations}
From assumption \eqref{assum positive}, we have
\bea \label{log C estimate}
\begin{cases}
 \llbracket{\ln C_i^\eps }\rrbracket \llbracket{C_i^\eps }\rrbracket=
(\ln C_{i,I}^\eps-\ln C_{i,E}^\eps)(C_{i,I}^\eps-C_{i,E}^\eps)=
\dfrac{1}{\xi}(C_{i,I}^\eps-C_{i,E}^\eps)^2 \geq \dfrac{1}{C_{u}}(C_{i,I}^\eps-C_{i,E}^\eps)^2  \\
 |\llbracket{\ln C_i^\eps }\rrbracket|=|\ln C_{i,I}^\eps-\ln C_{i,E}^\eps|=
|\dfrac{1}{\xi}(C_{i,I}^\eps-C_{i,E}^\eps)| \leq \dfrac{1}{C_{d}}|C_{i,I}^\eps-C_{i,E}^\eps|,
\end{cases}
\eea
where $\xi $ is between $C_{i,I}^\eps$ and $C_{i,E}^\eps$.
Then let $\varphi_s=\frac{\lambda_i}{z_i}C_{i,s}^\eps,~s=I,E$ in \eqref{weak form concentration two ion sub b}, we have
\bean
&&\eps \int_{\Gamma^\eps} \dfrac{\lambda_{i}}{z_i}P_m \pa_t\llbracket{\phi^\eps}\rrbracket \llbracket{C_i^\eps }\rrbracket ~ dS \no
&=& -\dfrac{\lambda_i}{z_i} \int_{\Omega^\eps} \sigma^\eps D \nabla \phi^\eps \nabla C_i^\eps~dx
-\eps \dfrac{\lambda_i}{z_i} \int_{\Gamma^\eps} \big( \sum\limits_{k=1}^{3}G_k\llbracket{\phi^\eps}\rrbracket \llbracket{C_i^\eps }\rrbracket+\sum\limits_{k=1}^{3}\dfrac{G_k}{z_k}\llbracket{\ln C_k^\eps }\rrbracket \llbracket{C_i^\eps }\rrbracket +I_p^\eps \llbracket{C_i^\eps }\rrbracket \big) ~ dS
\eean
Combining the above equaiton,   \eqref{log C estimate}, \eqref{assum positive}  and  \eqref{prior middle01 sub a} yields 
\bea \label{prior middle C}
&&\sum\limits_{i=1}^{3} ( \pa_t \| C_i^\eps \|^2_{\Omega^\eps}+
  \| \nabla C_i^\eps \|^2_{\Omega^\eps} +
\eps   \| \llbracket{C_i^\eps }\rrbracket \|^2_{\Gamma^\eps}) \no
&\leq&  A_1 \left( K_1 \| \nabla \phi^\eps  \|^2_{\Omega^\eps} +
\dfrac{1}{K_1} \sum\limits_{i=1}^{3}\| \nabla C_i^\eps  \|^2_{\Omega^\eps} \right) +
A_2 \eps \| \llbracket{\phi^\eps }\rrbracket \|^2_{\Gamma^\eps}+
A_3 \eps \sum\limits_{i=1}^{3}  \| \llbracket{C_i^\eps }\rrbracket \|^2_{\Gamma^\eps}
+A_4.
\eea
where we used the fact that $I_P^\eps$ is bounded. 

From assumption \eqref{assum positive 01}, we can multiply \eqref{prior middle01 sub b} with $\frac{2 A_1 K_1}{ D C_l }$ and add it to \eqref{prior middle C} to get
\bean
&&\sum\limits_{i=1}^{3} (\pa_t \| C_i^\eps \|^2_{\Omega^\eps}+
  \| \nabla C_i^\eps \|^2_{\Omega^\eps}+
\eps  \|  \llbracket{C_i^\eps }\rrbracket \|^2_{\Gamma^\eps})
 +
 A_1 K_1\| \nabla \phi^\eps \|^2_{\Omega^\eps} +
 A_2 \eps \| \llbracket{\phi^\eps }\rrbracket \|^2_{\Gamma^\eps}+
 A_3 \eps \pa_t\| \llbracket{\phi^\eps }\rrbracket \|^2_{\Gamma^\eps}   \no
&\leq& \dfrac{A_4}{K_1} \sum\limits_{i=1}^{3}\| \nabla C_i^\eps \|^2_{\Omega^\eps} +
A_5\eps \| \llbracket{\phi^\eps }\rrbracket \|^2_{\Gamma^\eps}+
A_6 \eps \sum\limits_{i=1}^{3}  \| \llbracket{C_i^\eps }\rrbracket \|^2_{\Gamma^\eps}+A_7 \no
&\leq & \dfrac{A_4}{K_1} \sum\limits_{i=1}^{3}\| \nabla C_i^\eps  \|^2_{\Omega^\eps} +
A_5\eps \| \llbracket{\phi^\eps }\rrbracket \|^2_{\Gamma^\eps}+
A_{8} \eps \sum\limits_{i=1}^{3}  \left(\eps^{-1} K_2 \| C_i^\eps \|^2_{\Omega^\eps}+
\dfrac{1}{K_2} \eps \|\nabla C_i^\eps \|^2_{\Omega^\eps} \right)+A_7,
\eean
where $K_1,K_2$ are two positive constants, and the trace theorem \cite{cioranescu1999introduction, Monsurro(2003)} has been used
\be  \label{boundary est02}
\left\|v_{s}\right\|_{L^{2}\left(\Gamma^{\epsilon}\right)} \leq A\left( \delta  \eps^{-\frac{1}{2}}\left\|v_{s}\right\|_{L^{2}\left(\Omega_{s}^\eps\right)}+
\dfrac{1}{\delta }\eps^{\frac{1}{2}}\left\|\nabla v_{s}\right\|_{L^{2}\left(\Omega_{s}^\eps\right)}\right),
\ee
for $\forall~ v_s \in H^1(\Omega_s^\eps) ,~s=I,E$. 

Choosing $K_1,K_2$  big enough, we have
\bea
&&\sum\limits_{i=1}^{3} ( \pa_t \| C_i^\eps \|^2_{\Omega^\eps}+
 \| \nabla C_i^\eps \|^2_{\Omega^\eps} +
\eps   \|  \llbracket{C_i^\eps }\rrbracket \|^2_{\Gamma^\eps})+
\| \nabla \phi^\eps \|^2_{\Omega^\eps} +
  \eps \| \llbracket{\phi^\eps }\rrbracket \|^2_{\Gamma^\eps}+
  \eps \pa_t\| \llbracket{\phi^\eps }\rrbracket \|^2_{\Gamma^\eps}   \no
&\leq& A_1 \bigg( \sum\limits_{i=1}^{3} \| C_i^\eps \|^2_{\Omega^\eps} + \eps \| \llbracket{\phi^\eps }\rrbracket \|_{\Gamma^\eps}^2\bigg)+A_2.
\eea
  Then by using Gronwall inequality, we obtain 
\bea \label{prior result01}
&&\sup_{t \in [0,T]} (\sum\limits_{i=1}^{3} \| C_i^\eps \|^2_{\Omega^\eps}+
\eps \| \llbracket{\phi^\eps }\rrbracket \|_{\Gamma^\eps}^2 )+
\sum\limits_{i=1}^{3}  \int_0^T \|\nabla C_i^\eps \|^2_{\Omega^\eps}~dt\no
&+& \int_0^T \|\nabla \phi^\eps \|^2_{\Omega^\eps}~dt +
\sum\limits_{i=1}^{3}  \int_0^T \eps \| \llbracket{C_i^\eps }\rrbracket \|_{\Gamma^\eps}^2~dt+
\int_0^T \eps \| \llbracket{\phi^\eps }\rrbracket \|_{\Gamma^\eps}^2~dt \leq A.
\eea
This leads to \eqref{th prior result sub a} and \eqref{th prior result sub b}. 

For the $L^2$ estimates of $\phi^\eps$, by assumption  \eqref{phi_s assump zero mean}, Lemma \ref{lemma poincare for phi} and \eqref{prior result01}, we obtain
\ben
\int_0^T \big( \| \phi_I^\eps \|^2_{\Omega_I^\eps}+\| \phi_E^\eps \|^2_{\Omega_E^\eps} \big)~dt
\leq A_1 \int_0^T \big( \|\nabla \phi_I^\eps \|^2_{\Omega_I^\eps}+\|\nabla \phi_E^\eps \|^2_{\Omega_E^\eps}+
\eps \| \llbracket{\phi^\eps }\rrbracket \|^2_{\Gamma^\eps} \big)~dt \leq A_2,
\een
\eqref{th prior result sub f} is obtained by combining this with \eqref{prior result01}. 

Using again the inequality \eqref{boundary est02}, \eqref{th prior result sub c} follows from \eqref{th prior result sub b} and \eqref{th prior result sub e} follows from \eqref{th prior result sub f}.

Lastly we show $\|\pa_t C^\eps_{i,s} \|_{L^2(0,T;H^{-1}(\Omega^\eps_s))}, i=1,2,3,s=I,E$ are bounded, we only consider $s=I$ since other cases can be discussed similarly. From \eqref{weak form concentration two ion sub b} we have
\begin{eqnarray}
    &&\dfrac{\eps \lambda_i}{z_i} P_m \int_{\Gamma^\eps}  \pa_t\llbracket{\phi^\eps}\rrbracket v_I~dS\nonumber\\
  &=&
  -\dfrac{ \lambda_i}{z_i} \int_{\Omega_I^\eps}   \sigma_I^\eps D \nabla \phi_I^\eps \nabla v_I~dx
  -\dfrac{\eps \lambda_i}{z_i}\int_{\Gamma^\eps}  \sum\limits_{k=1}^{3} \left(
    G_k \llbracket{\phi^\eps}\rrbracket v_I+\dfrac{G_k}{z_k}
\llbracket{\ln C_k^\eps }\rrbracket v_I +I_p^\eps v_I  \right)~dS .
\end{eqnarray}
 
Substituting the above equation into \eqref{weak form concentration two ion sub a} leads to
\bean
 &&\|\pa_t C^\eps_{1,I} \|_{H^{-1}(\Omega^\eps_I)} \\
 &\leq & A_1 \sup\limits_{\substack{v_1 \in H^1(\Omega^\eps_I),\\ \|v_1\|_{H^1(\Omega^\eps_I)}\leq 1} }
 \left(( \| \nabla C^\eps_{1,I}\|_{\Omega^\eps_{I}}+ \| \nabla \phi^\eps_{I}\|_{\Omega^\eps_{I}})  \| \nabla v_1\|_{\Omega^\eps_{I}} +
 \eps \| \llbracket{\phi^\eps }\rrbracket \|_{\Gamma_\eps} \|v_1\|_{\Gamma_\eps}\right.\\
&&\left. +
 \eps \sum\limits_{i=1}^{3} \| \llbracket{C_i^\eps }\rrbracket \|_{\Gamma_\eps} \|v_1\|_{\Gamma_\eps}+
 \sqrt{\eps} \|v_1\|_{\Gamma_\eps} \right) \\
 &\leq & A_2 \big( \| \nabla C^\eps_{1,I}\|_{\Omega^\eps_{I}} +\| \nabla \phi^\eps_{I}\|_{\Omega^\eps_{I}}
 +\sqrt{\eps} \| \llbracket{\phi^\eps }\rrbracket \|_{\Gamma_\eps}+
 \sqrt{\eps}  \| \llbracket{C_1^\eps }\rrbracket \|_{\Gamma_\eps} \big)+A_3,
\eean
where we used $\|v_1\|_{\Gamma_\eps} = O(\eps^{-1/2})$. 
Finally,  combining above equation with \eqref{prior result01}  yields 
\bean
&&\|\pa_t C^\eps_{1,I} \|^2_{L^2(0,T;H^{-1}(\Omega^\eps_I))} \\
&\leq & A_1 \int_0^T \left( \| \nabla C^\eps_{1,I}\|^2_{\Omega^\eps_{I}}+
 \| \nabla \phi^\eps_{I}\|_{\Omega^\eps_{I}}+
 \eps  \| \llbracket{\phi^\eps }\rrbracket \|^2_{\Gamma_\eps} +
 \eps\sum\limits_{i=1}^{3}   \| \llbracket{C_i^\eps }\rrbracket \|^2_{\Gamma_\eps}  \right) ~dt+A_2 \\
 &\leq & A_3.
\eean
\hspace*{\fill} $\square$

%
\section{Homogenization}  \label{section homo 2}
In this section,  the homogenization theories for the EN model \eqref{elec neu nondimen two ions} and \eqref{potential nodimen two ion} are developed with different connectivity.  The homogenization process for ``connected-disconnected'' case is shown in subsection \ref{subsec con-discon}, and then the ``connected-connected'' case is shown in subsection \ref{subsec con-con}.

\subsection{``Connected-disconnected'' case} \label{subsec con-discon}
In order to use compactness results from unfolding operators and two-scale theory, we extend the functions $C_{i,s}^\eps,\phi_s^\eps,i=1,2,3,s=I,E$ in a suitable way to the whole domain $\Omega$. Since $\Omega_E^\eps$ is connected and has a Lipschitz boundary, by \cite[Theorem 2.2]{hopker2014note}, there exists a linear and bounded extension operator  $\mathcal{L}^\eps:~L^{2}\left((0, T); H^{1}\left(\Omega_{E}^{\epsilon}\right)\right) \rightarrow L^{2}\left((0, T); H^{1}(\Omega)\right)$, we simply denote
 \be
 \widetilde{C}_{i,E}^\eps=\mathcal{L}^\eps C_{i,E}^\eps,~ \widetilde{\phi}_{E}^\eps=\mathcal{L}^\eps \phi_{E}^\eps.
 \ee
 The a priori estimates \eqref{th prior result sub a}, \eqref{th prior result sub b} and \eqref{th prior result sub f} remain valid for the extensions $\widetilde{C}_{i,E}^\eps$ and $\widetilde{\phi}_{E}^\eps$. Such extensions can not be applied to the functions in $\Omega_I^\eps$ since it is not connected \cite{Monsurro(2003)}.   Zero extensions is used instead. The zero extension of a function defined on $\Omega_s^\eps$ for $s=I,E$ will be denoted by $\overline{C}_{i,I}^\eps$ and $\overline{\phi}_I^\eps$.
It is obvious that the time derivative of $\overline{C}_{i,I}^\eps$  satisfies
\bea
&\left\langle\partial_{t} \overline{C}_{i,s}^\eps, \psi\right\rangle_{\Omega}=\left\langle\partial_{t} {C}_{i,s}^\eps, \psi \right\rangle_{\Omega_{s}^{\epsilon}} \text { for all } \psi \in H^{1}(\Omega) \text { and a.e. } t \in(0, T),\\
&\left\|\partial_{t} \overline{C}_{i,s}^\eps\right\|_{L^{2}\left((0, T); H^{-1}(\Omega)\right)} \leq
\left\|\partial_{t} {C}_{i,s}^\eps\right\|_{L^{2}((0, T); H^{-1}(\Omega_{s}^{\epsilon}))} \leq A,
\eea
where $\langle \cdot, \cdot  \rangle_G$ represent the duality paring $\langle \cdot, \cdot  \rangle_{H^{-1}(G)\times H^1(G)}$ for an arbitrary measurable set $G$. Since unfolding operators and two-scale convergence method are used in our proof, we give some definitions and properties for them in the Appendix. First we give a lemma.
\begin{lemma} \label{lem three multi}
Suppose $a_\eps,b_\eps,c_\eps \in L^2(\omega)$ be three series such that
\ben
a_\eps \rightarrow a ~~\text{strongly in}~L^2(\omega),~~
b_\eps \rightarrow b ~~\text{strongly in}~L^2(\omega),~~
c_\eps \rightharpoonup c ~~\text{weakly in}~L^2(\omega),
\een
and $\|a_\eps\|_{L^\infty(\omega)}\leq A,~\|a\|_{L^\infty(\omega)}\leq A,~\|b_\eps\|_{L^\infty(\omega)}\leq A,~\|b\|_{L^\infty(\omega)}\leq A$ for some open subset $\omega \in \mathbb{R}^d$, then the following convergence result holds:
\bea  \label{assertion conv}
\int_\omega  a_\eps b_\eps c_\eps ~dx \rightarrow \int_\omega a b c~dx .
\eea
\textit{Proof}: We have
\bean
\int_\omega  a_\eps b_\eps c_\eps~dx
= \int_\omega  a b c ~dx + \int_\omega a b (c_\eps-c)~dx +\int_\omega a (b_\eps-b) c_\eps~dx
+\int_\omega (a_\eps-a) b_\eps c_\eps~dx ,
\eean
it is obvious that the limits of the last three terms on the right hand side are zero. For example, for the last term on the right hand side, we have
\ben
\left| \int_\omega (a_\eps-a) b_\eps c_\eps~dx \right|
\leq A \int_\omega |(a_\eps-a) c_\eps| ~dx  \rightarrow 0,
\een
the assertion is proved.
\end{lemma}
The following theorem is the convergence and homogenization result for ``connected-disconnected'' case:

\begin{theorem} \label{theorem con-discon}
Suppose the intracellular region $\Omega_I^\eps$ is disconnected and assumption \eqref{phi_s assump zero mean} holds, let $C_{i,s}^\eps,\phi_s^\eps,~i=1,2,3,s=I,E$ be the solutions of \eqref{elec neu nondimen two ions} and \eqref{potential nodimen two ion}, then there exist subsequences of $C_{i,s}^\eps,\phi_s^\eps$, still denoted as $C_{i,s}^\eps,\phi_s^\eps$, and $C_{i,E} \in L^2(0,T;H^1(\Omega)),~\widehat{C}_{i,E} \in L^2((0,T)\times \Omega;H^1_{per}(Y))$, $C_{i,I} \in L^2((0,T)\times\Omega),~\widehat{C}_{i,I} \in L^2((0,T)\times \Omega;H^1(Y))$ and $\phi_E \in L^2(0,T;H^1(\Omega)),~\widehat{\phi}_E \in~L^2((0,T)\times \Omega;H^1_{per}(Y)),~ \phi_I \in L^2((0,T)\times \Omega),~\widehat{\phi}_I \in~L^2((0,T)\times \Omega;H^1(Y))$, such that, for $i=1,2,3$
\begin{subequations} \label{C_E conv01}
\begin{align}
\mathcal{T}^{\eps}(\widetilde{C}_{i,E}^\eps) \rightarrow C_{i,E}  ~~& \text{strongly in} ~L^2((0,T)\times \Omega\times Y),\label{C_E conv01 sub a}\\
\mathcal{T}^{\eps}(\nabla \widetilde{C}_{i,E}^\eps) \rightharpoonup \nabla_x C_{i,E} +\nabla_y \widehat{C}_{i,E}   ~~&
\text{weakly in} ~L^2((0,T)\times \Omega\times Y),\label{C_E conv01 sub b}\\
\mathcal{T}_b^{\eps}(\widetilde{C}_{i,E}^\eps) \rightarrow C_{i,E}  ~~& \text{strongly in} ~L^2((0,T)\times \Omega\times \Gamma),\label{C_E conv01 sub c}\\
\pa_t \overline{C}_{i,E}^\eps  \rightharpoonup |Y_E|\pa_t C_{i,E} ~~& \text{weakly in} ~L^2(0,T; H^{-1}(\Omega) ).\label{C_E conv01 sub d}
\end{align}
\end{subequations}
\vspace{-20pt}
\begin{subequations} \label{C_I conv01}
\begin{align}
\mathcal{T}^{\eps}(\overline{C}_{i,I}^\eps) \rightarrow \beta_I C_{i,I} ~~ & \text{strongly in} ~L^2((0,T)\times \Omega\times Y),\label{C_I conv01 sub a}\\
\mathcal{T}^{\eps}(\nabla\overline{ C_{i,I}^\eps}) \rightharpoonup  \beta_I \nabla_y \widehat{C}_{i,I}  ~~&
\text{weakly in} ~L^2((0,T)\times \Omega\times Y),\label{C_I conv01 sub b}\\
\mathcal{T}_b^{\eps}(\overline{C}_{i,I}^\eps) \rightarrow C_{i,I}  ~~& \text{strongly in} ~L^2((0,T)\times \Omega\times \Gamma),\label{C_I conv01 sub c}\\
\pa_t \overline{C}_{i,I}^\eps  \rightharpoonup |Y_I|\pa_t C_{i,I} ~~& \text{weakly in} ~L^2(0,T; H^{-1}(\Omega) ).\label{C_I conv01 sub d}
\end{align}
\end{subequations}
\vspace{-20pt}
\begin{subequations} \label{phi conv01}
\begin{align}
\mathcal{T}^{\eps}( \nabla \widetilde{\phi}_E^\eps ) \rightharpoonup \nabla_x \phi_E+\nabla_y \widehat{\phi}_E ~~& \text{weakly in}~L^2((0,T)\times \Omega\times Y),\label{phi conv01 sub a}\\
\mathcal{T}_b^{\eps}(\widetilde{\phi}_E^\eps ) \rightharpoonup \phi_E
~~& \text{weakly in}~L^2((0,T)\times \Omega\times \Gamma),\label{phi conv01 sub b}\\
\mathcal{T}^{\eps}( \overline{ \nabla \phi_I^\eps} ) \rightharpoonup \beta_I \nabla_y \widehat{\phi}_I
~~& \text{weakly in}~L^2((0,T)\times \Omega\times Y),\label{phi conv01 sub c}\\
\mathcal{T}_b^{\eps}({\phi}_I^\eps ) \rightharpoonup \phi_I
~~& \text{weakly in}~L^2((0,T)\times \Omega\times \Gamma).\label{phi conv01 sub d}
\end{align}
\end{subequations}
And $C_{i,s},\phi_s,~i=1,2,3,s=I,E$ are the solutions of the following equations
\bea \label{Th homo equ concen}
\begin{cases}
|Y_E| \pa_t C_{i,E} - \nabla \cdot (D_E^\star \nabla C_{i,E}+z_i C_{i,E} D_E^\star \nabla \phi_E)\\
=|\Gamma|\bigg( \dfrac{G_i}{z_i}(\phi_I-\phi_E)+
\dfrac{G_i}{z^2_i}\ln \dfrac{C_{i,I}}{C_{i,E}}+\dfrac{P_i^0}{z_i}
+\dfrac{\lambda_i}{z_i} P_m \pa_t (\phi_I-\phi_E)\bigg) & in~(0,T)\times \Omega,\\
-(D_E^\star \nabla C_{i,E}+z_i C_{i,E} D_E^\star \nabla \phi_E) \cdot \vec{\nu}=0 & on~(0,T)\times \pa \Omega, \\
|Y_I| \pa_t C_{i,I} \\
= -|\Gamma|\bigg( \dfrac{G_i}{z_i}(\phi_I-\phi_E)+
\dfrac{G_i}{z^2_i}\ln \dfrac{C_{i,I}}{C_{i,E}}+\dfrac{P_i^0}{z_i}
+\dfrac{\lambda_i}{z_i} P_m \pa_t (\phi_I-\phi_E)\bigg) & in~(0,T)\times \Omega, \\
C_{i,I}|_{t=0}=C_{i,I}^0,~~C_{i,E}|_{t=0}=C_{i,E}^0 & in~\Omega,
\end{cases}
\eea
\bea  \label{Th homo equ phi}
\begin{cases}
-\nabla \cdot \big( \sigma_E D_E^\star \nabla \phi_E \big)=0 & in ~(0,T) \times\Omega,\\
-\sigma_E D_E^\star \nabla \phi_E \cdot \vec{\nu} =0  & on ~(0,T) \times\pa \Omega,\\
P_m \pa_t (\phi_I-\phi_E)=- \big( \sum\limits_{i=1}^3 G_i(\phi_I-\phi_E)+\sum\limits_{i=1}^3 \dfrac{G_i}{z_i}\ln \dfrac{C_{i,I}}{C_{i,E}} +I_p^0 \big)
 &in ~(0,T) \times \Omega,\\
(\phi_I-\phi_E)|_{t=0}=\phi^0 &in~\Omega,
\end{cases}
\eea
where $\sigma_E=\sum\limits_{i=1}^3 z_i^2 C_{i,E}$, $D_E^\star = \int_{Y_E} D(I+\nabla_y \chi_E)~dy$ and $\chi_E=(\chi_E^1,\chi_E^2,\cdots,\chi_E^d)$ is the solution of
\eqref{asist equ E}, $P_1^0=3 I_p^0,P_2^0=-2 I_p^0, P_3^0=0$, and $I_p^0$ has the same expression as in \eqref{I_p express} but the concentrations are replaced by the corresponding limits.
And $C_{i,s}$ are electro-neutral, i.e. $\sum\limits_{i=1}^3 z_i C_{i,s}=0$ for $s=I,E$.
\end{theorem}
\textit{Proof}:
Convergence \eqref{C_E conv01}, \eqref{C_I conv01} and \eqref{phi conv01} are the results of a priori estimate Theorem \ref{Th a prior est} and relavent convergence results in \cite{Gahn(2016)}. 
Note that $C_E,~C_I,~\phi_E,~\phi_I$ are independent of $y$. This property will be used in the following derivation process.

 First, in \eqref{weak form concentration two ion sub b}, let
\be  \label{test01}
\varphi_I=\eps u(x)v_1(x/\eps)w(t),
\ee
where $u \in C_0^\infty(\Omega),~v_1 \in C_{per}^\infty(Y_I),~w \in C_0^\infty(0,T)$.
Then integrating \eqref{weak form concentration two ion sub b} with respect to time gives
\bea \label{mid04}
0&=&2 \int_0^T \int_{\Omega_I^\eps} \sigma_I^\eps D \nabla \phi_I^\eps (\eps v_1 w \nabla u+ u w \nabla_y v_1)~dx dt+\int_0^T \int_{\Gamma_\eps} \eps^2 \sum\limits_{i=1}^3 G_i\llbracket{\phi^\eps }\rrbracket u v_1 w ~dS dt \no
&&+
  \int_0^T\!\!\! \int_{\Gamma_\eps} \eps^2  \sum\limits_{i=1}^3 \dfrac{G_i}{z_i}\llbracket{\ln C_i^\eps  }\rrbracket
 u v_1 w ~dS dt
  -  \int_0^T \!\!\!\int_{\Gamma_\eps} \eps^2 P_m \llbracket{\phi^\eps }\rrbracket u v_1 \pa_t w ~dS dt\nonumber\\
  &&+\int_0^T \!\!\!\int_{\Gamma_\eps} \eps^2 I_p u v_1 w ~dS dt \no
& =&M_1+M_2+M_3-M_4+M_5.
\eea
 For $M_1$, using the property of unfolding operator Lemma \ref{lemma unfolding properties 01}, we deduce
\bea \label{M_1}
M_1
=   \int_0^T \int_{\Omega\times Y_I} \mathcal{T}^{\eps}\left( \overline{\sigma}_I^\eps \right) D~
 \mathcal{T}^{\eps}\left(\overline{ \nabla \phi_I^\eps } \right) ~
 (\eps v_1 w \mathcal{T}^{\eps}\left( \nabla u \right) +
 \mathcal{T}^{\eps}\left( u \right)w \nabla_y v_1)~ dx dy dt.
\eea
By \eqref{C_I conv01 sub a} and the relation \eqref{sigma_s relation}: $\sigma_I^\eps= \sum\limits_{i=1}^3 z_i^2 C_{i,I}^\eps$, it follows that
\bea \label{prove limit simga_I conv}
\mathcal{T}^{\eps}\left( \overline{\sigma}_I^\eps \right) \rightarrow  \beta_I \sigma_I~~~
 \text{strongly in}~ L^2((0,T)\times \Omega \times Y_I),
\eea
where $\sigma_I= \sum\limits_{i=1}^3 z_i^2 C_{i,I}$.
From Lemma \ref{lemma unfolding properties 03}, it is obvious that, as $\eps \rightarrow 0$
\bea \label{prove limit unfold test_I conv}
(\eps v_1 w \mathcal{T}^{\eps}\left( \nabla u \right) +
 \mathcal{T}^{\eps}\left( u \right)w \nabla_y v_1) \rightarrow u w \nabla_y v_1
 ~~~ \text{strongly in}~ L^2((0,T)\times \Omega \times Y_I).
\eea
From \eqref{prove limit unfold test_I conv}, \eqref{prove limit simga_I conv}, \eqref{phi conv01 sub c}, and Lemma \ref{lem three multi}, passing to the limit for \eqref{M_1} as $\eps \rightarrow 0$, while noticing that $C_0^\infty ((0,T)\times \Omega) \bigotimes C_{per}^\infty(Y_I) $ is dense in $L^2((0,T) \times\Omega;H^1(Y_I)) $, we deduce
\bean
M_1 \rightarrow   \int_0^T \int_{\Omega\times Y_I} \sigma_I D \nabla_y \widehat{\phi}_I \nabla_y \Phi_I ~ dx dy dt.
\eean
for $\forall~\Phi_I \in L^2((0,T) \times\Omega;H^1(Y_I))$.
For $M_2$, from Lemma \ref{lemma unfolding properties 01} we have
\bean
M_2=\int_0^T \int_{\Omega \times \Gamma} \sum\limits_{i=1}^{3} G_i
\big( \mathcal{T}_b^{\eps}\left( \phi_I^\eps  \right)-
\mathcal{T}_b^{\eps}( \widetilde{\phi}_E^\eps  ) \big) \eps \mathcal{T}_b^{\eps}\left( u\right)v_1 w ~ dx dS_y dt,
\eean
Obviously,
\ben
\eps \mathcal{T}_b^{\eps}\left( u\right)v_1 w~\rightarrow~0~~~\text{strongly in} ~L^2((0,T)\times \Omega\times \Gamma),
 \een
by \eqref{phi conv01 sub b} and \eqref{phi conv01 sub d} it follows that $M_2 \rightarrow 0$. For $M_4$, the result is similar, i.e. $M_4 \rightarrow 0$.

Then, for $M_3$, noticing that $\mathcal{T}_b^\eps(\ln C_{i,s}^\eps)=\ln \mathcal{T}_b^\eps(C_{i,s}^\eps)$, from Lemma \ref{lemma unfolding properties 01} it follows
\ben
M_3=\int_0^T \int_{\Omega\times\Gamma} \sum\limits_{i=1}^{3} (\ln \mathcal{T}_b^\eps(\overline{C}_{i,I}^\eps)-
\ln \mathcal{T}_b^\eps(\widetilde{C}_{i,E}^\eps)) \eps \mathcal{T}_b^\eps(u)v_1 w  ~ dx dS_y dt.
\een

  The results in \cite[Corrolary 15]{Gahn(2016)} yield that $\ln(\mathcal{T}_b^\eps(\widetilde{C}_{i,E}^\eps))$ and $\mathcal{T}_b^\eps(I_p^\eps)$ converges weakly in $L^2((0,T)\times \Omega\times \Gamma)$ to $\ln(C_{i,E})$ and $I_p^0$ respectively. Then we have   $\eps \mathcal{T}_b^\eps(u)v_1 w$ converges to zero strongly in $L^2((0,T)\times \Omega\times \Gamma)$. So the limits of $M_3$ and $M_5$ are zero as $\eps \rightarrow 0$.

Summarizing the above convergence results for $M_1,M_2,M_3,M_4,M_5$, we deduce that the limit equation for \eqref{mid04} as $\eps \rightarrow 0$ is
\be \label{con-discon pf mid01}
\int_0^T  \int_{\Omega\times Y_I} \sigma_I D \nabla_y  \widehat{\phi}_I \nabla_y \Phi_I ~ dx dy dt=0,
\ee
for $\forall ~\Phi_I \in L^2((0,T) \times\Omega;H^1(Y_I))$.
Let $\Phi_I=\widehat{\phi}_I$ and noticing $\sigma_I>0$, we obtain
\be \label{mid06}
\nabla_y \widehat{\phi}_I =0.
\ee
Similar results were derived in, for example \cite[Proposition 11]{Gahn(2016)} and \cite[Theorem 3.3]{Donato(2011)}.

Then let $\varphi_I=\varphi_1(t,x)$ in \eqref{weak form concentration two ion sub b}, where $\varphi_1 \in C_0^\infty((0,T) \times \Omega)$.  Using similar arguments as for $M_1,M_2,M_3,M_4,M_5$ while noticing \eqref{mid06}, we obtain, as $\eps \rightarrow 0$,
\bean 
\int_0^T \int_{\Omega \times \Gamma} \bigg[  \sum\limits_{i=1}^3 G_i(\phi_I-\phi_E)
+ \sum\limits_{i=1}^3 \dfrac{G_i}{z_i}(\ln C_{i,I}-\ln C_{i,E})
+I_p^0 + P_m \pa_t (\phi_I-\phi_E) \bigg] \varphi_1 ~ dx dS_y dt=0.
\eean
Since $\phi_I,~\phi_E$ and $C_{i,I},~C_{i,E}$ are independent of $y$, so
\bea  \label{phi homo01}
 \sum\limits_{i=1}^3 G_i(\phi_I-\phi_E)
+ \sum\limits_{i=1}^3 \dfrac{G_i}{z_i}(\ln C_{i,I}-\ln C_{i,E})
+I_p^0 + P_m \pa_t (\phi_I-\phi_E)=0.
\eea

If we let $\varphi_I=u(x)w_1(t)$ in \eqref{weak form concentration two ion sub b} and pass to the limit as $\eps \rightarrow 0$, where $w_1 \in C^\infty(0,T), w_1(T)=0$, then the initial condition for $\phi_I-\phi_E$ can be derived as
\bean
(\phi_I-\phi_E) |_{t=0}=\phi^0.
\eean

Then in \eqref{weak form concentration two ion sub b}, let $\varphi_E=\eps u(x)v_2(x/\eps)w(t)$, where $v_2 \in C_{per}^\infty(Y_E)$ and $u,w$ are the same as in \eqref{test01}. Using similar discussion as for \eqref{mid04}, from \eqref{phi conv01 sub a} and \eqref{assertion conv} we have, after passing to the limit for \eqref{weak form concentration two ion sub b} as $\eps \rightarrow 0$,
\be \label{mid05}
2\int_0^T \int_{\Omega \times Y_E} \sigma_E D (\nabla_x \phi_E+\nabla_y \widehat{\phi}_E) \nabla_y \Phi_E ~ dx dy dt=0,
\ee
for $\forall~ \Phi_E \in L^2((0,T) \times\Omega;H_{per}^1(Y_E))$. 

Above equation yields the following representation for $\widehat{\phi}_E$
\ben
\widehat{\phi}_E= \nabla_x \phi_E~\cdot \chi_E,
\een
where $\chi_E=(\chi_E^1,\chi_E^2,\cdots,\chi_E^d)$, satisfying the following auxiliary problem 
\bea  \label{asist equ E}
\begin{cases}
-\nabla_y\cdot \big(D\nabla_y(y_j+\chi_E^j) \big)=0~~~~~~in ~Y_E, \\
-D\nabla_y(y_j+\chi_E^j) \cdot \vec{n}_1=0~~~~~~~~~on~\Gamma,\\
\chi_E^j~ is~ Y-periodic,~~~~\int_{Y_E} \chi_E^j=0.
\end{cases}
\eea

If we  let $\varphi_E=\varphi_1(t,x)$ in \eqref{weak form concentration two ion sub b}, from \eqref{phi homo01} and \eqref{assertion conv}, using density property, it yields  the limit for \eqref{weak form concentration two ion sub b} as $\eps \rightarrow 0$  
\ben
2 \int_0^T \int_{\Omega} \sigma_E D_E^\star \nabla \phi_E \nabla \Phi_1 ~ dx dt=0,
\een
for $\forall~ \Phi_1 \in L^2((0,T) ; H^1(\Omega))$, where $D_E^\star = \int_{Y_E} D(I+\nabla_y \chi_E)~dy$. The strong form of the above equation is
\bea  \label{phi homo02}
\begin{cases}
-\nabla \cdot \big( \sigma_E D_E^\star \nabla \phi_E \big)=0 & in ~(0,T) \times \Omega,\\
-\sigma_E D_E^\star \nabla \phi_E \cdot \vec{\nu} =0 &on ~(0,T) \times\pa \Omega.
\end{cases}
\eea

 Now consider the equations for $C_i^\eps$.
Let $v_I=\eps u(x)v_1(x/\eps)w(t)$ in \eqref{weak form concentration two ion sub a}, where $u,v_1,w$ is the same as in \eqref{test01}. We deal with the interface integral similarly as for \eqref{mid04}. From \eqref{C_I conv01 sub b}, \eqref{assertion conv}, \eqref{mid06} and density property, passing to the limit for \eqref{weak form concentration two ion sub a} as $\eps \rightarrow 0$ leads to
\be \label{mid 106}
\int_0^T \int_{\Omega\times Y_I}  D \nabla_y  \widehat{C}_{i,I} \nabla_y \Phi_I ~ dx dy dt=0,
\ee
for $\forall ~\Phi_I \in L^2((0,T) \times\Omega;H^1(Y_I))$.
Let $\Phi_I=\widehat{C}_{i,I}$ in the above equation, we have
\be \label{mid grad C_I=0}
\nabla_y  \widehat{C}_{i,I}=0.
\ee

If let $v_I=\varphi_1(t,x)$ in \eqref{weak form concentration two ion sub a}, using \eqref{C_I conv01 sub d}, \eqref{assertion conv}, \eqref{mid06}, \eqref{mid grad C_I=0} and passing to the limit for \eqref{weak form concentration two ion sub a} as $\eps \rightarrow 0$ yields
\bean 
\hspace{-25pt}
&& |Y_I| \int_0^T \int_{\Omega} \pa_t C_{i,I} u w_1 ~ dx dt\no
\hspace{-25pt}
 &=& -  \int_0^T \int_{\Omega\times \Gamma} \bigg( \dfrac{G_i}{z_i}(\phi_I-\phi_E)
+\dfrac{G_i}{z^2_i}(\ln C_{i,I} -\ln C_{i,E})+\dfrac{P_i^0}{z_i}
+\dfrac{\lambda_i}{z_i} P_m \pa_t (\phi_I-\phi_E) \bigg)u w_1 ~ dx dS_y dt.
\eean

The integrand in the above equation is independent of $y$, so it is equivalent to
\be  \label{C homo02}
|Y_I| \pa_t C_{i,I} = -|\Gamma|\bigg( \dfrac{G_i}{z_i}(\phi_I-\phi_E)+
\dfrac{G_i}{z^2_i}(\ln C_{i,I} -\ln C_{i,E})+\dfrac{P_i^0}{z_i}
+\dfrac{\lambda_i}{z_i} P_m \pa_t (\phi_I-\phi_E) \bigg). 
\ee

To derive the initial conditions for $C_{i,I}$, if we let $v_I=u(x) w_1(t)$ in \eqref{weak form concentration two ion sub a} and pass to the limit as $\eps \rightarrow 0$, it is easy to show that $C_{i,I}|_{t=0}=C_{i,I}^0$.

For $C_{i,E}^\eps$, let $v_E=\eps u(x)v_2(x/\eps)w(t)$ in \eqref{weak form concentration two ion sub b}. From \eqref{assertion conv}, passing to the limit as $\eps \rightarrow 0$ leads to
\begin{eqnarray}
  \int_0^T \int_{\Omega \times Y_E} \left( D(\nabla C_{i,E}+\nabla_y \widehat{C}_{i,E})+z_i C_{i,E} D (\nabla_x \phi_E+\nabla_y \widehat{\phi}_E)  \right)\nabla_y \Phi_E~ dx dy dt 
\end{eqnarray}
 
for $\forall ~\Phi_E \in L^2((0,T) \times\Omega;H_{per}^1(Y_E))$.

Thanks to \eqref{mid05}, the second term on the left hand side is zero, and  the following representation holds,
\ben
\widehat{C}_{i,E}= \nabla_x C_{i,E}~\cdot \chi_E.
\een

Lastly, if let $v_E=\varphi_1(t,x)$ in \eqref{weak form concentration two ion sub a} and pass to the limit as $\eps \rightarrow 0$, we have by density property that
\bean
&&\int_0^T  \int_\Omega \left( |Y_E| \pa_t C_{i,E} \Phi_1+ D_E^\star \nabla C_{i,E} \nabla \Phi_1+
z_i   C_{i,E} D_E^\star \nabla \phi_E \nabla \Phi_1 \right) ~ dx dt\no
 -&& \hspace{-5pt} \int_0^T  \int_{\Omega \times \Gamma}\bigg( \dfrac{G_i}{z_i}(\phi_I-\phi_E)+
\dfrac{G_i}{z^2_i}(\ln C_{i,I} -\ln C_{i,E})+\dfrac{P_i^0}{z_i}
+\dfrac{\lambda_i}{z_i} P_m \pa_t (\phi_I-\phi_E)\bigg) \Phi_1 ~ dx dS_y dt=0,
\eean
for $\forall~ \Phi_1 \in L^2((0,T) ; H^1(\Omega))$.

The derivation of initial condition for $C_{i,E}$ is the same as for $C_{i,I}$, so we omit it. The strong form for the above equation is
\bea  \label{C homo01}
\begin{cases}
|Y_E| \pa_t C_{i,E} - \nabla \cdot (D_E^\star \nabla C_{i,E}+ z_i C_{i,E} D_E^\star \nabla \phi_E)=|\Gamma|\bigg( \dfrac{G_i}{z_i}(\phi_I-\phi_E)+\\

\dfrac{G_i}{z^2_i}(\ln C_{i,I} -\ln C_{i,E}) +\dfrac{P_i^0}{z_i}+
\dfrac{\lambda_i}{z_i} P_m \pa_t (\phi_I-\phi_E)\bigg) &  in~(0,T)\times \Omega,\\
-(D_E^\star \nabla C_{i,E}+z_i C_{i,E} D_E^\star \nabla \phi_E) \cdot \vec{\nu}=0 & on~(0,T)\times \pa \Omega,\\
C_{i,E}|_{t=0}=C_{i,E}^0   &  in ~\Omega.
\end{cases}
\eea
\par Now we have derived all the homogenization equations, i.e. \eqref{C homo01}, \eqref{C homo02}, \eqref{phi homo01}, \eqref{phi homo02}, which are \eqref{Th homo equ concen} and \eqref{Th homo equ phi}. The electro-neutrality condition for $C_{i,s}$ is obvious since the initial value of $C_{i,s}$ is electro-neutral. The theorem is proved. \hspace*{\fill}   $\square$

\subsection{``Connected-connected'' case}  \label{subsec con-con}
Now we consider the ``connected-connected'' case. In this case, the status of $\Omega_I^\eps$ and $\Omega_E^\eps$ are ``equivalent'', which means the homogenized equations for both regions are in the same form. Noticing that $\Omega_I^\eps$ also reaches the boundary of $\pa \Omega$, we need to add boundary conditions for $C_I^\eps,~\phi_I^\eps$ on $\pa \Omega$, i.e. $J_{i,I}\cdot \vec{\nu}=0,~\text{and}~D \sigma_{I}^\eps  \nabla \phi_I^\eps  \cdot \vec{\nu}=0  ,~~   on~\pa \Omega$.

Since now both $\Omega_I^\eps$ and $\Omega_E^\eps$ are connected, by \cite[Theorem 2.2]{hopker2014note}, there exist linear and bounded extension operators $\mathcal{L}_s^\eps:~L^{2}\left((0, T), H^{1}\left(\Omega_{s}^{\epsilon}\right)\right) \rightarrow L^{2}\left((0, T), H^{1}(\Omega)\right), s=I,E$. Simply denote
 \be
 \widetilde{C}_{i,s}^\eps=\mathcal{L}_s^\eps C_{i,s}^\eps,~ \widetilde{\phi}_{s}^\eps=\mathcal{L}_s^\eps \phi_{s}^\eps,~~~s=I,E.
 \ee
The priori estimates \eqref{th prior result sub a}, \eqref{th prior result sub b} and \eqref{th prior result sub f} remain valid for the extensions $\widetilde{C}_{i,s}^\eps$ and $\widetilde{\phi}_{s}^\eps$.
The convergence and homogenization results for ``connected-connected'' case is summarized  in the  following theorem.

\begin{theorem}   \label{two ions converge th 01}
Suppose $\Omega_I^\eps$ is connected and assumption \eqref{phi_s assump zero mean} holds, let $C_{i,s}^\eps,\phi_s^\eps,~i=1,2,3,s=I,E$ be solutions of \eqref{elec neu nondimen two ions} and
\eqref{potential nodimen two ion}, then there exist subsequences of $C_{i,s}^\eps,\phi_s^\eps$, still denoted as $C_{i,s}^\eps,\phi_s^\eps$, and $C_{i,s},\phi_s \in L^2(0,T;H^1(\Omega)),~\widehat{C}_{i,s},\widehat{\phi}_s \in L^2((0,T)\times \Omega;H^1_{per}(Y))$, such that for $i=1,2,3,s=I,E$,
\begin{subequations} \label{C conv01 con-con}
\begin{align}
\mathcal{T}^{\eps}(\widetilde{C}_{i,s}^\eps) \rightarrow C_{i,s} ~~ & \text{strongly in} ~L^2((0,T)\times \Omega\times Y),\label{C conv01 con-con sub a}\\
\mathcal{T}^{\eps}(\nabla \widetilde{C}_{i,s}^\eps) \rightharpoonup \nabla_x C_{i,s} +\nabla_y \widehat{C}_{i,s}  ~~&
\text{weakly in} ~L^2((0,T)\times \Omega\times Y),\label{C conv01 con-con sub b}\\
\mathcal{T}_b^{\eps}(\widetilde{C}_{i,s}^\eps) \rightarrow C_{i,s} ~~ & \text{strongly in} ~L^2((0,T)\times \Omega\times \Gamma),\label{C conv01 con-con sub c}\\
\pa_t \overline{C}_{i,s}^\eps  \rightharpoonup |Y_s|\pa_t C_{i,s} ~~& \text{weakly in} ~L^2(0,T; H^{-1}(\Omega) ).\label{C conv01 con-con sub d}
\end{align}
\end{subequations}
\begin{subequations} \label{phi conv01 con-con}
\begin{align}
\mathcal{T}^{\eps}( \nabla \widetilde{\phi}_s^\eps ) \rightharpoonup \nabla_x \phi_s+\nabla_y \widehat{\phi}_s~~ & \text{weakly in}~L^2((0,T)\times \Omega\times Y),\label{phi conv01 con-con sub a}\\
\mathcal{T}_b^{\eps}(\widetilde{\phi}_s^\eps ) \rightharpoonup \phi_s
~~& \text{weakly in}~L^2((0,T)\times \Omega\times \Gamma).\label{phi conv01 con-con sub b}
\end{align}
\end{subequations}
 $C_{i,s},\phi_s,~i=1,2,3,s=I,E$ are the solutions of the following equations
\bea   \label{Th01 homo equ concen}
\begin{cases}
|Y_s| \pa_t C_{i,s} - \nabla \cdot (D_s^\star \nabla C_{i,s}+z_i C_{i,s} D_s^\star \nabla \phi_s)\\
=\pm |\Gamma|\bigg( \dfrac{G_i}{z_i}(\phi_I-\phi_E)  +
\dfrac{G_i}{z^2_i}\ln \dfrac{C_{i,I}}{C_{i,E}}+\dfrac{P_i^0}{z_i}
+\dfrac{\lambda_i}{z_i} P_m \pa_t (\phi_I-\phi_E)\bigg) & in~(0,T)\times \Omega,\\
(D_s^\star \nabla C_{i,s}+z_i C_{i,s} D_s^\star \nabla \phi_s) \cdot \vec{\nu}=0 & on~(0,T)\times \pa \Omega, \\
C_{i,s}|_{t=0}=C_{i,s}^0  & in~\Omega,
\end{cases}
\eea
\bea   \label{Th01 homo equ phi}
\begin{cases}
-\nabla \cdot \big( \sigma_s D_s^\star \nabla \phi_s \big) \\
=\pm \bigg(\sum\limits_{i=1}^{3} G_i(\phi_I-\phi_E)+\sum\limits_{i=1}^{3}\dfrac{G_i}{z_i}\ln \dfrac{C_{i,I}}{C_{i,E}}+I_p^0 +P_m \pa_t (\phi_I-\phi_E)  \bigg) & in ~(0,T) \times \Omega,\\
\sigma_s D_s^\star \nabla \phi_s \cdot \vec{\nu} =0  & on ~(0,T) \times\pa \Omega,\\
(\phi_I-\phi_E)|_{t=0}=\phi^0 &in~\Omega,
\end{cases}
\eea
where ``$\pm$'' takes ``$+$'' for $s=E$ and ``$-$'' for $s=I$. $D_s^\star = \int_{Y_s} D(I+\nabla_y \chi_s)~dy$, $\sigma_s=\sum\limits_{i=1}^3 z_i^2 C_{i,s}$ and $\chi_s=(\chi_s^1,\chi_s^2,\cdots,\chi_s^d),s=I,E$ are the solutions of
\eqref{asist equ E} and \eqref{two ions intra cell pro con-con}, $P_i^0,I_p^0$ have the same meaning as in theorem \ref{theorem con-discon}. And $C_{i,s}$ are electro-neutral, i.e. $\sum\limits_{i=1}^3 z_i C_{i,s}=0$ for $s=I,E$.
\end{theorem}
\textit{Proof}: The proof of the convergence results \eqref{C conv01 con-con} and \eqref{phi conv01 con-con} are the same as the proof of the convergence results for $\widetilde{C}_{i,E}^\eps,\widetilde{\phi}_E^\eps$ in Theorem \ref{theorem con-discon}, so we omit it. Since the status of $\Omega_I^\eps$ and $\Omega_E^\eps$ are ``equivalent'' in ``connected-connected'' case, we can see that the convergence results are the same for intracellular functions and extracellular functions. The homogenized equations \eqref{Th01 homo equ concen} and \eqref{Th01 homo equ phi} are different from the results in Theorem \ref{theorem con-discon}, but the proof are quite similar, so we only give a simplified proof and point out the differences. As the derivation of homogenized equations in $\Omega_I^\eps$ and $\Omega_E^\eps$ are the same, we only consider equations in $\Omega_I^\eps$.

 First consider $\phi_I^\eps$. Let $\varphi_I=\eps u(x)v_1(x/\eps)w(t)$ in \eqref{weak form concentration two ion sub b}, where $u \in C_0^\infty(\Omega),~v_1 \in C_{per}^\infty(Y_I),~w \in C_0^\infty(0,T)$. By convergence result \eqref{C conv01 con-con} and \eqref{phi conv01 con-con}, we can use similar argument as for \eqref{mid04} to get that, for $\forall~ \Phi_I \in L^2((0,T) \times\Omega;H_{per}^1(Y_I))$,
\be \label{mid 105}
\int_0^T \int_{\Omega \times Y_I} \sigma_I D (\nabla_x \phi_I+\nabla_y \widehat{\phi}_I) \nabla_y \Phi_I~ dx dy dt=0.
\ee
This result is different from \eqref{con-discon pf mid01} in last subsection, which is due to convergence result \eqref{phi conv01 con-con sub a} for $s=I$. 

From \eqref{mid 105} we have the following representation
\ben
\widehat{\phi}_I= \nabla_x \phi_I~\cdot \chi_I,
\een
where $\chi_I=(\chi_I^1,\chi_I^2,\cdots,\chi_I^d)$, satisfying
\bea   \label{two ions intra cell pro con-con}
\begin{cases}
-\nabla_y\cdot \big(D\nabla_y(y_j+\chi_I^j) \big)=0~~~~~~in ~Y_I ,\\
-D\nabla_y(y_j+\chi_I^j) \cdot \vec{n}_1=0~~~~~~~~~on~\Gamma.\\
\chi_I^j~ is~ Y-periodic,~~~~\int_{Y_I} \chi_I^j=0
\end{cases}
\eea

Then let $\varphi_I=\varphi_1(t,x)$ in \eqref{weak form concentration two ion sub b}, where $\varphi_1 \in C_0^\infty((0,T) \times \Omega)$, from convergence result \eqref{C conv01 con-con} and \eqref{phi conv01 con-con}, by density property, we can pass to the limit for \eqref{weak form concentration two ion sub b} as $\eps \rightarrow 0$ to get
\bean
 &&\int_0^T \int_{\Omega} \sigma_I D_I^\star \nabla \phi_I \nabla \Phi_1~ dx dt \\
    &-&\int_0^T \int_{\Omega \times \Gamma} \bigg[ \sum\limits_{i=1}^3 G_i(\phi_I-\phi_E)
+   \sum\limits_{i=1}^{3} \dfrac{G_i}{z_i}\ln \dfrac{C_{i,I}}{C_{i,E}}
+I_p^0 + P_m \pa_t (\phi_I-\phi_E) \bigg] \Phi_1 ~dx dS_y dt=0,
\eean
for $\forall \Phi_1 \in L^2((0,T) ; H^1(\Omega))$, where $D_I^\star = \int_{Y_I} D(I+\nabla_y \chi_I)~dy$. The corresponding strong form is
\bea   \label{phi homo02 con-con}
\begin{cases}
-\nabla \cdot \big( \sigma_I D_I^\star \nabla \phi_I \big) \\
=- \bigg(\sum\limits_{i=1}^3 G_i(\phi_I-\phi_E)+\sum\limits_{i=1}^{3}\dfrac{G_i}{z_i}\ln \dfrac{C_{i,I}}{C_{i,E}}+I_p^0 +P_m \pa_t (\phi_I-\phi_E)  \bigg) & in ~(0,T) \times \Omega,\\
\sigma_I D_I^\star \nabla \phi_I \cdot \vec{\nu} =0  & on ~(0,T) \times\pa \Omega,
\end{cases}
\eea
which is \eqref{Th01 homo equ phi} for $s=I$.

Next we consider $C_{i,I}^\eps$. Let $v_I=\eps u(x)v_1(x/\eps)w(t)$ in \eqref{weak form concentration two ion sub a}, passing to the limit as $\eps \rightarrow 0$ leads to
\ben
\int_0^T \int_{\Omega \times Y_I} \left( D(\nabla C_{i,I}+\nabla_y \widehat{C}_{i,I}) +
z_i  C_{i,I} D (\nabla_x \phi_I+\nabla_y \widehat{\phi}_I) \right)  \nabla_y \Phi_I ~ dx dy dt=0,
\een
for $\forall ~\Phi_I \in L^2((0,T) \times\Omega;H_{per}^1(Y_I))$. This is different from \eqref{mid 106} in last subsection, which is due to the convergence result \eqref{phi conv01 con-con sub a} and \eqref{C conv01 con-con sub b} for $s=I$.

According to \eqref{mid 105}, the second term on the left hand side of the above equation is zero, so we have
\ben
\widehat{C}_{i,I}= \nabla_x C_{i,I}~\cdot \chi_I,
\een
where $\chi_I$ is the solution of \eqref{two ions intra cell pro con-con}. Then, let $v_I=\varphi_1(t,x)$ in \eqref{weak form concentration two ion sub a}, passing to the limit as $\eps \rightarrow 0$ yields
\bean
&&\int_0^T  \int_\Omega \left(|Y_I| \pa_t C_{i,I} u_1 w +  D_I^\star \nabla C_{i,I} \nabla u_1 w+
z_i  C_{i,I} D_I^\star \nabla \phi_I \nabla u_1 w \right)~ dx dt \no
 &+&\int_0^T \int_{\Omega \times \Gamma}\bigg( \dfrac{G_i}{z_i}(\phi_I-\phi_E)+
\dfrac{G_i}{z^2_i}\dfrac{\ln C_{i,I}}{\ln C_{i,E}} +\dfrac{P_i^0}{z_i} +\dfrac{\lambda_i}{z_i} P_m \pa_t (\phi_I-\phi_E)\bigg) u_1 w ~ dx dS_y dt=0.
\eean
The corresponding strong form is
\bea  \label{C homo01 con-con}
\begin{cases}
|Y_I| \pa_t C_{i,I} - \nabla \cdot (D_I^\star \nabla C_{i,I}+z_i C_{i,I} D_I^\star \nabla \phi_I)\\
=-|\Gamma|\bigg( \dfrac{G_i}{z_i}(\phi_I-\phi_E)  +
\dfrac{G_i}{z^2_i}\dfrac{\ln C_{i,I}}{\ln C_{i,E}}+\dfrac{P_i^0}{z_i} +
\dfrac{\lambda_i}{z_i} P_m \pa_t (\phi_I-\phi_E)\bigg) & in~(0,T)\times \Omega,\\
-(D_I^\star \nabla C_{i,I}+z_i C_{i,I} D_I^\star \nabla \phi_I) \cdot \vec{\nu}=0 & on~(0,T)\times \pa \Omega.
\end{cases}
\eea
which is \eqref{Th01 homo equ concen} for $s=I$. The initial condition for $C_{i,I}$ can be obtained by similar arguments as in the proof of Theorem \ref{theorem con-discon}.  Summarizing the above results, we have \eqref{Th01 homo equ concen}, \eqref{Th01 homo equ phi}. The electro-neutrality condition for $C_{i,s}$ is a result of the electro-neutrality of $C_{i,s}^0$. The theorem is proved. \hspace*{\fill}   $\square$

\section{Conclusion}
In this paper, a micro-scale EN bi-domain ion transport model is first proposed to model the ion transportation on tissues. 
On the membrane, we considered the passive currents induced by ion channels, active currents induced by pumps and currents induced by the capacitance property of membrane.  
Then the homogenization  theory  for the nonlinear coupled system is derived rigorously by using unfolding operator and two-scale convergence method.  Different connectivity conditions for intracellular region $\Omega_I^\eps$ leads to different homogenized equations. When both the intracelluar region and extracellular region are connected to form a syncytium, repectively, the obtained system is a diffusion-convection-reaction system. While, when the intracellulars are disconnected, the macroscale effective equation of intracellular region is only  a reaction equation. In this case, the intracellular could indirectly communicate through the connected  extracellular space. Our macroscale model could be used to study some diseases induced by ion micro-circulation disorder, like spreading depression \cite{chang2013mathematical} and lens problems \cite{zhu2019bidomain}.
\appendix
\section{ Two-scale convergence and unfolding operator} 
We present in this Appendix some homogenization theory and results about two-scale convergence and unfolding operator. The two-scale convergence theory in periodic domains was first establish by Nguetseng \cite{nguetseng1989general} and further developed by Allaire \cite{Allaire(1992)}. Later two-scale convergence was extended to periodic interfaces \cite{Allaire(1996),Neuss-Radu(1996)}. First we give the definition of two-scale convergence with time ``t'' and strong two-scale convergence. Let $\omega$ be an open set in $\mathbb{R}^d$, and $Y=(0,1)^d$.
\begin{definition}
 Let $u_{\eps} \subset L^{2}((0, T) \times \omega),~u_{0} \in L^{2}((0, T) \times \omega \times Y)$, $u_{\eps}$ is said to converge (weakly) in the two-scale sense to $u_{0}$, if for $\forall \phi \in C\left([0, T] \times \overline{\omega}, C_{p e r}(\overline{Y})\right)$ the following relation holds:
\be \label{two scale defi001}
\lim _{\eps \rightarrow 0} \int_{0}^{T} \int_{\omega} u_{\eps}(t, x) \phi\left(t, x, \frac{x}{\eps}\right) ~d x d t=\int_{0}^{T} \int_{\omega} u_{0}(t, x, y) \phi(t, x, y) ~d y d x d t.
\ee
 If, in addition to \eqref{two scale defi001}, the following relation holds:
\be \label{two scale defi002}
\lim _{\eps \rightarrow 0}\left\|u_{\eps}\right\|_{L^{2}((0, T) \times \omega)}=\left\|u_{0}\right\|_{L^{2}((0, T) \times \omega \times Y)},
\ee
then $u_{\eps}$ is said to converge strongly in the two-scale sense to $u_{0}$.
\end{definition}
The following compactness result is an extension of the stationary case in \cite{Allaire(1992)} to the case with time ``t''. The time variable can be treated as an independent parameter and all the proof remain valid.
\begin{lemma} \label{two scale th01}
(i) Every bounded sequence in $L^{2}((0, T) \times \omega)$ has a two-scale convergent subsequence.\\
(ii) Let $u_{\eps}$ be a bounded sequence in $L^{2}\left((0, T), H^{1}(\omega)\right)$. Then there exists
$u_0~\in~L^{2}$ \\ $\left((0, T), H^{1}(\omega)\right)$, $u_{1} \in L^{2}\left((0, T) \times \omega, H_{p e r}^{1}(Y) \right)$ and a subsequence $u_{\eps}$ still denoted by $u_{\eps}$, such that
\ben
u_{\eps} \rightarrow u_{0},~~~\nabla u_{\eps} \rightarrow \nabla_{x} u_{0}+\nabla_{y} u_{1},~~~~\text{in the two-scale sense}.
\een
\end{lemma}
 
Next we introduce unfolding operator. Unfolding operator was first constructed in \cite{Arbogast(1990)} and further studied in detail in \cite{Cioranescu(2002),Cioranescu(2008)}. There are some equivalence results between two-scale convergence and unfolding operator.
 Let $\Gamma$ be a (n-1)-dimensional Lipschitz manifold, compactly included in $Y$, and
$\Gamma_\eps=\bigcup\limits_{k\in \mathbb{Z}^n} \eps(\Gamma+k),~ \Gamma_\eps~\subset~\omega$. Denote $[x]$ as the largest integer not greater than $x$, and let $ \{x\}:=x-[x] $, then we have
$[x] \in \mathbb{Z}^{n},~\{x\} \in \overline{Y}$. Let $Y^{*} \subset Y$ and $\omega_\eps$ be the internal of
$\omega \cap \bigcup\limits_{k \in \mathbb{Z}^{n}} \eps(k+\overline{Y^{*}})$.
\begin{definition} \label{definition unfolding op}
(i) Let $u_{\eps} \in L^{2}((0, T) \times \omega)$, unfolding operator
$\mathcal{T}^{\eps}: L^{2}((0, T) \times \omega) \rightarrow L^{2}((0, T) \times \omega \times Y)$ is defined by
\ben
\left(\mathcal{T}^{\eps} u_{\eps}\right)(t, x, y):=u_{\eps}\left(t, \eps\left(\left[\frac{x}{\eps}\right]+y\right)\right).
\een
(ii) Similarly, unfolding operator in $\omega_\eps$ is defined by
$\mathcal{T}^{\eps}: L^{2}\left((0, T) \times \omega_{\eps}\right) \rightarrow
L^{2}((0, T) $ \\ $\times \omega \times Y^{*})$, and unfolding operator on $\Gamma_\eps$ is defined by
$\mathcal{T}^{\eps}_{b}: L^{2}((0, T) \times\left.\Gamma_{\eps}\right) \rightarrow L^{2}((0, T) \times \omega \times \Gamma) \text { on } \Gamma_{\eps}$.
\end{definition}
From the above definition it is clear that $\mathcal{T}^{\eps}$ and $\mathcal{T}^{\eps}_{b}$ are bounded linear operators. Some properties of unfolding operator are given in the following lemma
\begin{lemma} \label{lemma unfolding properties 01}
\cite{Cioranescu(2002),Cioranescu(2008),Cioranescu(2006)}
(i) Let $u_{\eps}, v_{\eps} \in L^{2}\left((0, T) \times \omega_{\eps}\right)$, then we have
\ben
\left(\mathcal{T}^{\eps} u_{\eps}, \mathcal{T}^{\eps} v_{\eps}\right)_{L^{2}((0, T) \times \Omega \times Y^\star)}=\left(u_{\eps}, v_{\eps}\right)_{L^{2}\left((0, T) \times \omega_{\eps}\right)},
\een
and for $u_{\eps} \in L^{2}\left((0, T), H^{1}\left(\omega_{\eps}\right)\right)$ we have
 $\nabla_{y} \mathcal{T}^{\eps} u_{\eps}=\eps \mathcal{T}^{\eps}\left(\nabla_{x} u_{\eps}\right)$.\\
(ii) Let $u_{\eps}, v_{\eps} \in L^{2}\left((0, T) \times \Gamma_{\eps}\right)$, then we have
\ben
\begin{aligned}\left(\mathcal{T}^{\eps}_{b} u_{\eps}, \mathcal{T}^{\eps}_{b} v_{\eps}\right)_{L^{2}((0, T) \times \Omega \times \Gamma)} &=\eps\left(u_{\eps}, v_{\eps}\right)_{L^{2}\left((0, T) \times \Gamma_{\eps}\right)} , \\\left\|\mathcal{T}^{\eps}_{b} u_{\eps}\right\|_{L^{2}((0, T) \times \Omega \times \Gamma)} &=\sqrt{\eps}\left\|u_{\eps}\right\|_{L^{2}\left((0, T) \times \Gamma_{\eps}\right)}, \end{aligned}
\een
\end{lemma}
where $(\cdot,\cdot)_{L^2(G)}$ is the inner product of $L^2(G)$ for an arbitrary measurable set $G$. The relation between unfolding operator and two-scale convergence is given below (see \cite[Proposition 4.6]{Bourgeat(1997)} and \cite[Proposition 4.7]{Neuss-Radu(2007)}).
\begin{lemma} \label{lemma unfolding properties 02}
For bounded sequences in $L^{2}((0, T) \times \omega)$, the following two statements are equivalent\\
(i) $u_{\eps} \rightarrow u_{0}$ weakly/strongly in the two-scale sense,     \\
(ii) $\mathcal{T}^{\eps} u_{\eps} \rightarrow u_{0}$ weakly/strongly in $L^{2}((0, T) \times \Omega \times Y)$.\\
The same result holds for $u_\eps$ in $L^{2}\left((0, T) \times \Gamma_{\eps}\right)$ with
$\sqrt{\eps}\left\|u_{\eps}\right\|_{L^{2}\left((0, T) \times \Gamma_{\eps}\right)} \leq A$ and bounded unfolding operator $\mathcal{T}^{\eps}_{b}$.
\end{lemma}
For the unfolding operator in $\omega_\eps$, we have the following convergence result \cite[Proposition 2.6]{yang2014periodic}.
\begin{lemma} \label{lemma unfolding properties 03}
Let $\mathcal{T}^{\eps}: L^{2}\left((0, T) \times \omega_{\eps}\right) \rightarrow
L^{2}\left((0, T) \times \omega \times Y^{*}\right)$ and $u \in L^{2}\left((0, T)\times \omega\right)$, then we have
\ben
\mathcal{T}^{\varepsilon}(u) \rightarrow u ~~~\text {strongly in } L^{2}\left((0, T) \times \omega \times Y^{*} \right).
\een
\end{lemma}
We also state a lemma which is used in the proof of the a priori estimate \eqref{th prior result}.
\begin{lemma} \label{lemma poincare for phi}
Let $\phi_s^\eps,s=I,E$ be the solutions of \eqref{weak form concentration two ion}, if $\int_{\Omega_E^\eps} \phi_E^\eps=0$, then we have
\be  \label{poincare for phi}
\int_0^T \big( \| \phi_I^\eps \|^2_{\Omega_I^\eps}+\| \phi_E^\eps \|^2_{\Omega_E^\eps} \big)~dt
\leq A \int_0^T \big( \|\nabla \phi_I^\eps \|^2_{\Omega_I^\eps}+\|\nabla \phi_E^\eps \|^2_{\Omega_E^\eps}+
\eps \| \llbracket{\phi^\eps }\rrbracket \|^2_{\Gamma^\eps} \big)~dt
\ee
\end{lemma}
\textit{Proof}: Using poincar\'{e} inequality for functions with zero mean \cite{poincare1890equations,poincare1894equations}, we can
prove \eqref{poincare for phi} along the same line as the proof of \cite[Lemma 2.8]{Monsurro(2003)}.\hspace*{\fill} $\square$

\noindent  \textbf{Acknowledgment} This work was partially supported by the NSFC (\textcolor{red}{grant numbers}   11971342, 12071190), and NSERC (CA).
The first author is grateful to the financial support by China Scholarship Council (CSC) during his visit to York University in Toronto, where this work is done. \vspace{40pt}


%


{\footnotesize
\bibliographystyle{plain}
\bibliography{referenceii}
}

\end{document}